\documentclass[11pt]{article}

\input{amssym.def}
\input{amssym}
\usepackage{eufrak}
\usepackage{color}

\newtheorem{theorem}{Theorem}

\newtheorem{corollary}[theorem]{Corollary}
\newtheorem{definition}[theorem]{Definition}

\newtheorem{proposition}[theorem]{Proposition}
\newtheorem{remark}[theorem]{Remark}

\def\P{{\cal{P}}}
\def\pa{{\partial}}

\def\RRR{{\cal{R}}}

\def\om{\omega}

\def\Y{\mathbf{Y}}
\def\YR{\mathbf{Y}(R)}

\def\qq{q^{-1}}

\def\lhq{{\cal L}(\h, q)}

\def\h{\hbar}

\def\lr{{\cal L}(R)}

\def\LL{\tilde{L}}
\def\tlr{\tilde{\cal L}(R)}

\def\agl{\widehat{gl(V_R)}}

\def\De{\Delta}

\def\GG{{\frak{G}}}

\def\de{\delta}

\def\gg{\mbox{$\frak g$}}
\def\ot{\otimes}

\def\K{{\Bbb K}}

\def\C{{\Bbb C}}
\def\R{{\Bbb R}}

\def\Sym{{\rm Sym\, }}

\def\tr{{\rm tr}}
\def\End{{\rm End}}
\def\Ren{R^{\End(V)}}
\def\vv{V^{\otimes 2}}

\def\tLL{\tilde L}

\def\lhq{\ifmmode {\cal L}(q,\hbar)\else ${\cal L}(q,\hbar)$\fi}

\def\lqh{\ifmmode {\cal L}(q,\hbar)\else ${\cal L}(q,\hbar)$\fi}

\def\Tr{{\rm Tr}}

\def\he{\hat{e}}

\def\OT{\overline{T}}

\def\be{\begin{equation}}
\def\ee{\end{equation}}

\begin{document}

\makeatletter
\renewcommand{\theequation}{{\thesection}.{\arabic{equation}}}
\@addtoreset{equation}{section} \makeatother

\title{From Reflection Equation Algebra to Braided Yangians}

\author{\rule{0pt}{7mm} Dimitri
Gurevich\thanks{gurevich@ihes.fr}\\
{\small\it LAMAV, Universit\'e de Valenciennes,
59313 Valenciennes, France}\\
{\small \it and}\\
{\small \it   Interdisciplinary Scientific Center J.-V. Poncelet}\\
\rule{0pt}{7mm} Pavel Saponov\thanks{Pavel.Saponov@ihep.ru}\\
{\small\it
National Research University Higher School of Economics,}\\
{\small\it 20 Myasnitskaya Ulitsa, Moscow 101000, Russian Federation}\\
{\small \it and}\\
{\small \it
Institute for High Energy Physics, NRC "Kurchatov Institute"}\\
{\small \it Protvino 142281, Russian Federation}}

\maketitle

\begin{abstract}
In general, quantum matrix algebras are associated with a couple of compatible braidings. A particular example of such an algebra is
the so-called Reflection Equation algebra. In this paper we analyse its specific properties, which distinguish it from other quantum matrix
algebras (in first turn, from the RTT one). Thus, we exhibit a specific form of  the Cayley-Hamilton identity for its generating matrix, which in a  limit turns into the
Cayley-Hamilton identity  for the generating matrix of the enveloping algebra $U(gl(m))$. Also, we consider some specific properties of
the braided Yangians, recently introduced by the authors. In particular, we establish an analog of the Cayley-Hamilton identity for the
generating matrix of such a braided Yangian. Besides, by passing to a limit of the braided Yangian, we get a Lie algebra similar to that entering
the construction of the rational Gaudin model. In its enveloping algebra  we construct a Bethe subalgebra by the
method due to D.Talalaev.
\end{abstract}

{\bf AMS Mathematics Subject Classification, 2010:} 81R50


\section{Introduction}

Let $V$ be a vector space, $\dim V=N$, and $R:\vv\to \vv$ be a braiding, i.e. a solution of the braid relation (also called the quantum Yang-Baxter equation)
$$
R_{12} R_{23} R_{12} =R_{23} R_{12} R_{23}, \qquad R_{12}=R\ot I,  \quad R_{23}=I\ot R.
$$
Hereafter $I$ stands for the identity operator or its matrix. The above relation is written in the space $V^{\otimes 3}$ and
the lower indices label the spaces  where a given operator acts.

The unital associative algebra generated by entries   of a matrix $L=\|l_i^j\|_{1\leq i,j\leq N}$,
subject to the following system
\be
R\, L_1\, R\,L_1-L_1\, R\, L_1\, R=0,\quad L_1=L\ot I,
\label{RE}
\ee
is called the {\em Reflection Equation} (RE) algebra, associated with a braiding $R$ and denoted $\lr$. Below, the matrix $L$ and all
similar matrices are called {\em  generating}.

The algebra $\lr$ is a particular case of the so-called quantum  matrix (QM) algebras. Any QM algebra is associated with a couple of
compatible braidings (see \cite{IOP} for more details).

Another well-known example of a QM algebra is the so-called RTT algebra, generated by entries of a matrix  $T=\|t_i^j\|_{1\leq i,j\leq N}$
subject to the system
\be
R\, T_1\, T_2-T_1\,  T_2\, R=0,\qquad   T_1=T\ot I,\quad T_2=I\ot T.
\label{RTT}
\ee
All algebras, we are dealing with, are assumed to be unital.

A braiding $R$ is called an {\em involutive symmetry}, if it meets the condition $R^2=I$, and  a {\em Hecke symmetry}, if it meets the
Hecke relation\footnote{From now on, the notation  $\K$ stands for the ground field, which is  $\C$ or $\R$.}
\be
(qI-R)(\qq I+R)=0,\qquad q\in \K,\quad q^2\not= 1.
\label{Hec}
\ee
If an involutive or Hecke symmetry  $R$ is a deformation of the flip $P$, then the both QM algebras are deformations of the commutative
algebra ${\rm Sym}(gl(N))$. This means that the dimensions of homogeneous components of each of these algebras are classical, i.e.
equal to those of corresponding components\footnote{If $R$ is a Hecke symmetry we should additionally require $q$ to be generic, that is
$q^n\not=1$ for any integer $n$.} in ${\rm Sym}(gl(N))$. Emphasize that similar algebras can be associated with any braiding $R$ but in general this
deformation property fails. Below, all symmetries  $R$ which are deformations of the usual flips and the corresponding objects will be
referred to as {\em deformation} ones.

The best known examples of deformation  Hecke symmetries are those coming from the Quantum Groups (QG) $U_q(sl(N))$. However,
we introduce all our algebras  without any QG, which plays merely the role of a symmetry group for them, provided the corresponding
Hecke symmetry $R$  comes from this QG. As another example of a deformation Hecke symmetry we mention the Cremmer-Gervais
$R$-matrices. However, in general the involutive and Hecke symmetries, we are dealing with, are not deformation either of the usual flips
or of the super-ones.

Also, we assume all symmetries to be skew-invertible (see the next section for the definition). Under this condition a braided
(or $R$-)analog of the trace can be defined. Note that this trace enters  all our constructions. In particular, it takes part in the
definition of quantum analogs of the symmetric polynomials in all algebras under consideration. These quantum symmetric polynomials
generate commutative subalgebras, called {\em characteristic}.

However, only  in the RE algebras these subalgebras are central (see \cite{IP}). Besides, the RE algebras possess many other
properties distinguishing them from other QM algebras. The main  purpose  of the present paper is to exhibit specific features of the RE
algebras and of the so-called {\it braided Yangians} \cite{GS1}, which are current (i.e. depending on parameters) algebras in a sense close to the RE ones.

Here, we mention two of these particular properties. First, if a Hecke symmetry $R=R(q)$ is a deformation of the usual flip $P$ (i.e.
$R(1)=P$), then the corresponding {\it modified} RE algebra (\ref{mRE}) is a deformation of the enveloping algebra $U(gl(N))$. It can be treated as the enveloping algebra
of a {\em braided analog} $gl(V_R)$ of a Lie algebra $gl(n)$.
If $R$ is a skew-invertible
Hecke symmetry of a general type, similar analogs of  the Lie algebra $gl(N)$ and its enveloping algebra can be also defined.

In this connection the following question arises: whether it is possible to define an affine version of the braided Lie algebras similar to
$\widehat{gl(N)}$? Below, we introduce such a braided analog  $\widehat{gl(V_R)}$. Note that, putting aside the affine QG
$U_q({\widehat \gg})$, there are known two approaches to define quantum generalizations of affine algebras: the RE algebras in the spirit
of \cite{RS}\footnote{They differ from our braided Yangians by the middle terms, which are also current $R$-matrices. Observe that there
are known many versions of such RE algebras.} and the double Yangians and their q-analogs as introduced in \cite{FJMR}. In our
subsequent publications we plan to study the centre of the enveloping algebra  $U(\widehat{gl(V_R)})$ in the frameworks of the Kac's approach and to compare all mentioned methods of defining quantum affine algebras.

The second particular property of the RE algebra is that its generating matrix $L$ satisfies a matrix polynomial identity  $Q(L)=0$ for
a polynomial $Q(t)$, called {\em characteristic}. Thus, we get a version of the Cayley-Hamilton identity.

As was shown in \cite{IOP}, such an identity exists for the generating matrices of other QM algebras. However, only in the RE algebra
this identity arises from the characteristic polynomial. Also,  in deformation cases by passing to the limit $q\to 1$ in  the modified form
of the RE algebra, we get the characteristic polynomial and the corresponding Cayley-Hamilton identity for the generating matrix of the enveloping algebra
$U(gl(N))$,  which are usually obtained via  the so-called Capelli determinant.

As noticed above, the   braided Yangians, recently introduced   in \cite{GS1}, are in a sense close to the RE algebras. They are associated with current  quantum $R$-matrices, constructed by means of the Yang-Baxterization of involutive and Hecke symmetries.
These braided Yangians constitute one of two classes which generalize   the Yangian $\Y(gl(N))$, introduced by V.Drinfeld \cite{D}.
The second class of Yangian-like algebras, also introduced in \cite{GS1}, consists of the so-called {\em Yangians of RTT type} which are more similar to the RTT algebras.

One of the main dissimilarities of the braided Yangians and these of RTT type arises from
their evaluation morphisms. For the braided Yangians the evaluation morphisms are  similar to these for the  Yangians  $\Y(gl(N))$, but  their  target algebras are the
RE algebras (modified or not) instead of $U(gl(N))$. Another particular property of the braided Yangians is that they admit the Cayley-Hamilton
identities for the generating matrices, which   are  also more similar to the classical ones. This is due to the fact that the analogs of the matrix powers
entering these identities are given by the usual matrix product of several copies of the generating matrix (but with a shifted parameter $u$).

Also, deformation braided Yangians, in particular those, associated with $R$-matrices (\ref{Rm1}) and called {\em braided $q$-Yangians}, admit a limit\footnote{In order to get this limit
we first change the basis in the Yangian, or in other words, we pass to  the shifted form of this Yangian.}
 as $q\to 1$. In this limit we get
 Lie algebras $\GG_{trig}$ similar to  $\GG$ entering construction of the rational Gaudin model. By using the method due to D.Talalaev, we construct  Bethe subalgebras in the
enveloping algebras $U(\GG_{trig})$.
Consequently, we get new Bethe subalgebras in the Lie algebras $gl(N)^{\oplus K}$. In a more detailed way the corresponding version of an integrable model
 will be considered elsewhere.

The paper is organized as follows. In the next section we recall some basic properties of braidings and symmetries. In section 3 we
describe the RE algebra and the corresponding braided Lie algebra $gl(V_R)$. Also, we  define its affinization. In Section 4 we consider
different forms of the characteristic polynomials for the generating matrices of the RE algebras. In Section 5 we introduce braided
Yangians and describe their specific properties. In the last section by passing to the $q=1$ limit in the braided $q$-Yangian,
we get the aforementioned current Lie algebra $\GG_{trig}$ and  find a Bethe subalgebra in its enveloping algebra.

\medskip

\noindent
{\bf Acknowledgement} D.G. is grateful to the Max Planck Institute for Mathematics (Bonn), where the paper was mainly written,  for
stimulating atmosphere during his  scientific visit. The work of P.S. has been funded by the Russian Academic Excellence Project '5-100' and was also partially supported by the RFBR grant 16-01-00562. The authors are also thankful to D.Talalaev for elucidating discussion.

\section{Braidings: definitions and properties}

The starting object of our approach is a {\em skew-invertible} braiding of one of two types, specified below. Recall that a braiding $R$ is
called skew-invertible if there exists an operators $\Psi:\vv\to \vv$ such that the following relation holds
\be
{\rm Tr}_2 R_{12} \Psi_{23}=P_{13}\quad\Leftrightarrow\quad R_{ij}^{kl} \Psi_{lp}^{jq}=\de_i^q \de_p^k,
\label{Psi}
\ee
where the symbol $\Tr_2$ means that the trace is applied in the second matrix space. Below a summation over repeated indexes
is always understood. Here we assume that a basis $\{x_i\}$ in the space $V$ is fixed and $\|R_{ij}^{kl}\|$ is the matrix of the operator
$R$ in the basis $\{x_i\ot x_j\}$:
$$
R(x_i\ot x_j)=R_{ij}^{kl}\,x_k\ot x_l.
$$

The condition (\ref{Psi}) enables us to extend $R$ up to a braiding
$$
R: \vv\to \vv,\quad  (V^*)^{\ot 2}\to (V^*)^{\ot 2},\quad V^*\ot V\to V\ot V^*,\quad V\ot V^*\to  V^*\ot V,
$$
such that there exists an $R$-invariant pairing $V\ot V^*\to \K$ (see \cite{GPS2}). This extended braiding implies a braiding
$$
\Ren: \,{\rm End}(V)^{\ot 2}\to {\rm End}(V)^{\ot 2},
$$
where we identify ${\rm End}(V)\cong V\ot  V^*$ since $V$ is a finite dimensional space.

Now, introduce two operators
\be
B={\rm Tr}_1 \Psi \quad \Leftrightarrow\quad B_i^j=\Psi_{ki}^{kj},\qquad\, C={\rm Tr}_2 \Psi\quad \Leftrightarrow\quad
C_i^j=\Psi_{ik}^{jk}.
\label{BC}
\ee
The definition of $\Psi$ and the Yang-Baxter equation for $R$ leads to the properties:
\be
{\rm Tr}_1\, B_1R_{12}=I_2,\qquad {\rm Tr}_2\, C_2 R_{12} =I_1.
\label{pro2}
\ee
\be
R_{12}B_1B_2 =B_1B_2R_{12},\qquad R_{12}C_1C_2 =C_1C_2R_{12},
\label{pro1}
\ee

Let $\{x^j\}$ be the {\it right} dual basis of the space $V^*$, i.e. $<x_i, x^j>=\de_i^j$. Then the $R$-invariant  pairing in the opposite order is
\be
<\,,\,>:\,V^*\ot V\to \K,\quad  <x^j, x_i>=B_i^j.
\label{par}
\ee

In the space $\End (V)$ we fix the following basis
$$
l_i^j:= x_i\ot x^j\in \End (V)
$$
and consider the map
\be
{\rm tr}_R:\,\End(V)\to \K,\quad l_i^j\mapsto \de_i^j,
\label{def:tr}
\ee
motivated by the pairing $V\ot V^*\to \K$. We call this map the $R$-trace.

Also, using the pairing (\ref{par}), we define the following product in the space  $\End(V)$
\be
\circ :\, \End(V)^{\ot 2}\to \,\End(V),\quad  l_i^j\circ l_k^l=B_k^j\, l_i^l.
\label{prod}
\ee
This product is $\Ren$-invariant in the following sense
\be
\Ren (I\ot \circ) = (\circ\ot I)\Ren_2  \Ren_1, \qquad \Ren (\circ \ot I)= (\circ \ot I)\Ren_1  \Ren_2,
\label{inv}
\ee
where all operators act on the space  $\End(V)^{\ot 3}$. Hereafter, for the sake of simplicity we write $R_k$ instead of $R_{k\,k+1}$.

Now, introduce the following pairing on the space $\End(V)$:
\be
\langle \,,\, \rangle :\, \End(V)^{\ot 2}\to \K,\quad \langle X,Y \rangle={\rm tr}_R (X\circ Y),\quad  \forall\, X,Y\in \End(V).
\label{pair}
\ee
Thus, on the generators $l_i^j$ we have $\langle l_i^j, l_k^l\rangle=B_k^j \de_i^l$.

Now, let  $M=\|m_i^j\|$ be an $N\times N$ matrix. Define its $R$-trace as
$$
{\rm Tr}_R M=\Tr(C M)= M_i^j\, C_j^i.
$$
This definition is motivated as follows. With a matrix $M$ we associate an element $M_i^j x_j\ot \tilde{x}^i\in\End(V)$,  where  $\tilde{x}^i$ is the {\it left}
dual basis of the space $V^*$, i.e. $<\tilde{x}^i, x_j>=\de_j^i$. The right $R$-invariant pairing of
$\tilde x^i$ and $x_j$ reads $<x_j, \tilde{x}^i>=C_j^i$ (see \cite{GPS2} for details). So, the $R$-trace of $M$ is just the result
of applying this pairing to $M_i^j x_j\ot \tilde{x}^i$.

As was shown in \cite{O}, this $R$-trace has an important property: for any $N\times N$ matrix $M$ the following holds
$$
\Tr_{R(2)} R_1 M_1 R_1^{-1}=\Tr_{R(2)} R_1^{-1} M_1 R_1=I_1\Tr_R M.
$$
From now on, we use the following notation $\Tr_{i_1\dots i_k}=\Tr_{i_1}\dots\Tr_{i_k}$ where  $i_1<\dots <i_k$ and the same for $R$-traces.

Now, we assume $R$ to be a Hecke symmetry and $q$ to be generic. The corresponding constructions and results for involutive symmetries can be obtained by passing
to the limit $q\to 1$.

Note that for any Hecke symmetry $R$ the symmetric and skew-symmetric algebras
$$
{\rm Sym}_R(V)=T(V)/\langle {\rm Im}(qI-R)\rangle,\qquad {\bigwedge}_R(V)=T(V)/\langle {\rm Im}(\qq I+R)\rangle
$$
can be introduced. Since they are graded, the corresponding Poincar\'e-Hilbert series
$$
P_+(t)= \sum_k \dim {\Sym}_R^{\!(k)}(V)t^k,\qquad  P_-(t)=\sum_k \dim  {\bigwedge}_R^{\!(k)}(V)t^k,
$$
are well defined. Here the index $(k)$ labels the $k$-th order homogeneous components. According to \cite{H} the
Poincar\'e-Hilbert series  $P_{\pm}(t)$ are rational functions.

Emphasize that the above homogeneous components can also be defined via the projectors of symmetrization $\P_+$ (called below symmetrizers) and skew-symmetrization
$\P_-$ (skew-sym\-met\-ri\-zers)
$$
\P_+^{(k)}:\, V^{ \ot k}\to  {\rm Sym}_R^{\!(k)}(V),\qquad  \P_-^{(k)}: V^{ \ot k}\to {\bigwedge}_R^{\!(k)}(V).
$$
The latter operators can be introduced by a recursive relation:
\be
\P^{(k)}_- =\frac{1}{k_q} \,\P^{(k-1)}_-\left(q^{k-1} I-(k-1)_q R_{k-1}\right)\P^{(k-1)}_-,\qquad
k_q=\frac{q^k-q^{-k}}{q-q^{-1}},
\label{skew}
\ee
where we put by definition $\P_-^{(1)} = I$ and  assume that the skew-symmetrizer $\P^{(k)}_-$ is always applied at the positions
$1,2,\dots,k$. Formula (\ref{skew}) was proved in \cite{G} in a little bit different normalization of the Hecke symmetries.

Let us assume the rational function  $P_-(t)$ to be noncancellable. Let $m$ (respectively $n$) be the degree of its
numerator (respectively,  denominator). The ordered couple $(m|n)$ is called the {\em bi-rank} of the
symmetry $R$.

If $R$ is a skew-invertible symmetry (involutive or Hecke) and its by-rank is $(m|n)$ then the operators $B$ and $C$ (\ref{BC})
have a few additional properties:
\be
BC = q^{2(n-m)}  I,\qquad \Tr\, B=\Tr\, C= q^{n-m} (m-n)_q.
\label{add-prop}
\ee

\begin{proposition}
\label{pr1}
If $R$ is a skew-invertible Hecke symmetry and its bi-rank is $(m|0)$, then \rm
\be
\Tr_{R(k+1\dots m)} \P_-^{(m)}=q^{-m(m-k)} \frac{k_q!(m-k)_q!}{m_q!} \P_-^{(k)},
\label{tr-asym}
\ee
\em where we use the notation $k_q!=1_q 2_q \dots k_q$ and the standard agreement $0_q! = 1$.
\end{proposition}


\noindent
{\bf Proof.} The proof of this claim is a direct consequence of recurrence (\ref{skew}) and properties of $R$-trace.
Indeed, in virtue of the condition on the bi-rank and (\ref{add-prop}) we have
$$
\Tr_R I=\Tr\, C = q^{-m}m_q,
$$
while formula (\ref{pro2}) means that $\Tr_{R(k)}R_{k-1} =I_{k-1}$. Now we can calculate a typical trace:
$$
\Tr_{R(k)}\left(q^{k-1} I-(k-1)_qR_{k-1}\right)=\left(q^{k-m-1} m_q - (k-1)_q \right)I_{k-1}=q^{-m}(m-k+1)_q I_{k-1},
$$
where at the last step we used the relation
$$
q^a b_q-q^b a_q=(b-a)_q.
$$
Thus, we have
$$
\Tr_{R(m)}\left( q^{m-1} I - (m-1)_qR_{m-1}\right)=q^{-m} I_{m-1}
$$
 and consequently,
$$
\Tr_{R(m)}\P^{(m)}_-=\frac{q^{-m} }{m_q} \,\P^{(m-1)}_-.
$$

Upon applying the $R$-trace once more, we get
$$
\Tr_{R(m-1, m)}\P^{(m)}_-=q^{-2m}\,\frac{2_q}{m_q(m-1)_q} \,\P^{(m-2)}_-= q^{-2m} \frac{2_q! (m-2)_q!}{m_q!}\,\P^{(m-2)}_-.
$$
Now, using the reasoning by recursion,  we arrive to formula  (\ref{tr-asym}).\hfill \rule{6.5pt}{6.6pt}

\section{Braided Lie algebras and their affinization}

Now, consider a unital associative algebra generated by matrix elements of $N\times N$ matrix $\LL= \|l_i^j\|$
which obey the system of quadratic-linear relations:
\be
R\, \LL_1\, R\,\LL_1-\LL_1\, R\, \LL_1\, R=R\,\LL_1-\LL_1\, R.
\label{mRE}
\ee
We call this algebra the {\em modified RE algebra} and denote it $\tlr$.

If $q^2\not=1$ the algebras $\lr$ (\ref{RE}) and $\tlr$ are isomorphic to each other. The isomorphism is realized by the following map
\be
\LL \mapsto L+\frac{1}{q-\qq} I.
\label{iso}
\ee

Due to this reason we treat the algebra $\tlr $ as a modified form of the algebra $\lr $. In \cite{GPS2} we have constructed a
representation category of the algebra $\tlr $ similar to that of the algebra $U(gl(N))$. The isomorphism (\ref{iso}) enables us to
convert any $\tlr$-module into a $\lr$-one.

Here, we want  to mention only three $\tlr $-modules. The first one is the basic space $V$. The corresponding {\it vector representation}
is defined by
$$
\rho_V(l_i^j) \triangleright x_k=B_k^j x_i.
$$
where the notation $\triangleright$ stands for the action of a linear operator.

The second $\tlr $-module, called {\em covector}, is defined in the dual space $V^*$ by the following action on basis elements
$$
\rho_{V^*}(l_i^j)\triangleright x^k=- x^l R^{kj}_{li}.
$$

The third module, called {\em adjoint}, is identified with $V\otimes V^*$. The action of the elements $l_i^j$ on this module is defined by means of the following coproduct
\be
\De(l_i^j)= l_i^j\ot 1 + 1\ot l_i^j-(q-\qq) \sum_k l_i^k\ot l_k^j.
\label{cop}
\ee
Onto the whole algebra $\tlr $ this coproduct is extended by means of the braiding $\Ren $. In this sense we speak about
a {\em braided bi-algebra structure} of the algebra $\tlr $. The reader is referred to \cite{GPS2} for details. Note, that the
coproduct (\ref{cop}) arises from the braided structure of the RE algebra discovered by Sh.Majid \cite{M}.

Another way to define the adjoint representation is based on a {\em braided analog} of the Lie bracket. It is defined as follows.
The system quadratic-linear relations (\ref{mRE}) on the generators of the algebra $\tlr $ can be rewritten as
$$
l_i^j\ot l_k^l-\RRR(l_i^j\ot l_k^l)=[l_i^j, l_k^l],
$$
where
$$
\RRR: \End(V)^{\ot 2}\to \End(V)^{\ot 2}\quad {\rm and}\quad  [\,,\,]:\End(V)^{\ot 2}\to \End(V)
$$
are two operators. The operator $\RRR$ is defined below (see (\ref{RRR})).

Emphasize that if $q=1$ the operators $\RRR$ and $\Ren$ coincide with each other but for a generic $q$ it is not so.

Then the adjoint representation of the algebra $\tlr$ can be defined as followed
\be
\rho_{\End(V)}  (l_i^j) \triangleright l_k^l = [l_i^j, l_k^l].
\label{adj}
\ee

\begin{proposition} \label{prop2} The action {\rm (\ref{adj})} defines  a representation of the algebra $\tlr$. \end{proposition}

In order to prove  this claim it suffices to show that the action (\ref{adj}) coincides with that discussed above.   It can be also shown by straightforward computations.

\begin{definition} The space $\End(V)$ equipped with the operators $\RRR$ and $[\,,\,]$ is called  braided   Lie algebra
and is denoted $gl(V_R)$. \end{definition}

Also, the algebra $\tlr$ plays the role of the enveloping algebra of the braided Lie algebra $gl(V_R)$ in virtue
of Proposition \ref{prop2}.

Besides the property formulated in Proposition \ref{prop2}, the braided Lie algebra $gl(V_R)$ has the following features.

1. Its bracket $[\,,\,]$ is skew-symmetric in the following sense:  $[\,,\,]{\cal P}=0$.  Here ${\cal P}:\,\End(V)^{\ot 2}\to \End(V)^{\ot 2}$ is
a symmetrizer, constructed in \cite{GPS2}. Note that in comparison with the above projectors $\P_{\pm}$ acting in tensor powers of the space
$V$, the symmetrizer $\P$  acts in $({\rm span}_\K(l_i^j))^{\ot 2}$. Such symmetrizers and analogical skew-symmetrizers were constructed in
\cite{GPS2} for the tensor powers 2 and 3 of the space ${\rm span}_\K(l_i^j)$.

2. This bracket is $\Ren$-invariant in the same sense as in (\ref{inv}). In the case related to the QG $U_q(sl(N))$, this bracket is also
covariant with respect to the action of this QG.

As for (\ref{adj}), we treat it as a braided analog of the Jacobi identity. Note that if $R$ is an involutive symmetry, then the corresponding
Jacobi identity can be written under the following form, similar to the classical (or super-)one:
\be
[\,,\,][\,,\,]_{12}(I+\Ren_1 \Ren_2+\Ren_2\Ren_1)=0.
\label{Jac}
\ee
Also, note that the right hand side of (\ref{mRE}) can be obtained by
applying the product $\circ$ (\ref{prod}) to its left hand side. This follows form the relation
$$
\circ L_1 R_{12} L_1= L_1\Tr_{(1)}B_1R_{12} = L_1I_2.
$$

Now, introduce the following useful notation. We  put\footnote{\label{fo}  Note that the QM algebras as introduced in \cite{IOP}, are
defined in a similar way, but with the help of the second braiding $F$, in a sense compatible with $R$:
$L_{\overline k}=F_{k-1} L_{\overline{k-1}}F_{k-1}^{-1}$.}
\be
L_{\overline 1 }=L_{1},\quad L_{\overline k}=R_{k-1} L_{\overline{k-1}}R_{k-1}^{-1},\quad k\ge 2.
\label{ove}
\ee
With the use of this notation, we can rewrite the system (\ref{RE}) in a form similar to (\ref{RTT}):
\be
R_1 L_{\overline 1}\, L_{\overline  2}-L_{\overline 1}\,L_{\overline 2} R_1=0.
\label{sys}
\ee
The adjoint action can also be written as
\be
L_{\overline 1}\triangleright R_1\, L_{\overline 1}=L_{\overline 1}-R_1^{-1}\, L_{\overline 1}\,
R_1 \quad\Leftrightarrow \quad L_{\overline 1} \triangleright  L_{\overline 2}= L_{\overline 1}\,R_1^{-1}-R_1^{-1}\, L_{\overline 1} =
L_{\overline 1}\,R_1-R_1\, L_{\overline 1},
\label{ac}
\ee
where in the last equality we use the following consequence of Hecke condition (\ref{Hec}):
$$
R^{-1}=R-(q-\qq)I.
$$

As for the operators $\RRR$ and $\Ren$, they can be respectively presented as
\be
\RRR(L_{\overline 1}\ot L_{\overline 2})=R^{-1}(L_{\overline 1}\ot L_{\overline 2})R,\qquad
 \Ren(L_{\overline 1}\ot L_{\overline 2})=L_{\overline 2}\ot L_{\overline 1}.
\label{RRR}
\ee
Below, we use similar notations for dealing with the so-called braided Yangians.

Now, consider the element
$$
\ell=\Tr_R L=\Tr(CL)\in \End(V).
$$
In order to stress the difference between the $R$-traces $\Tr_R$ and $\tr_R$ (see (\ref{def:tr})) note that $\tr_R L=I$.
We have
$$
\tr_R\ell=C^k_k=\frac{(m-n)_q}{q^{m-n}}.
$$
Let us suppose that $m\not= n$ and consequently $tr_R\ell\not=0$. Then the elements
$$
f_i^j=l_i^j-\de_i^j\,\frac{\ell}{\tr_R\ell}
$$
are well defined and traceless:  $\tr_R f_i^j=0$. This enables us to define a braided analog $sl(V_R)$ of the Lie algebra $sl(N)$.
The reader is referred to \cite{GPS2} for details.

Now, consider the affinization procedure of the braided Lie algebras $gl(V_R)$. For the algebras $sl(V_R)$ it can be done in a
similar manner. Following the classical pattern, we consider the algebra
$$
gl(V_R)[t, t^{-1}]= gl(V_R)\ot \K[t, t^{-1}].
$$
This algebra is generated by elements $l_i^j[a]:=l_i^j\ot t^a$, $a\in {\Bbb Z}$. The braided Lie bracket in it is also defined
according to the classical pattern:
$$
[X[a], Y[b]]:= [X,Y][a+b], \quad X, Y\in gl(V_R).
$$

To construct the central extension of  $gl(V_R)[t, t^{-1}]$ we introduce a vector space
$$
gl(V_R)[t, t^{-1}]\oplus \K\, c
$$
where $c$ is a new generator commuting with any elements of $gl(V_R)[t, t^{-1}]$. Besides, we extend the action of the operator
$\RRR $ in a natural way
$$
\RRR (X[a] \ot Y[b])= \RRR (X\ot Y)[b][a], \quad \RRR (X[a]\ot c)=c \ot X[a],\quad \RRR (c \ot X[a])=X[a]\ot c.
$$
Here, the notation $\RRR (X\ot Y)[b][a]$ means that we attribute the label $b$ to the first factor and that $a$ to the second one.

Now, define the {\em affine braided Lie algebra}  $\agl $ by introducing the following bracket
$$
[X[a], c]=0,\quad  [X[a], Y[b]]=[X,Y][a+b]+ \om(X[a], Y[b])\, c,
$$
where
$$
\om(X[a], Y[b]):=a \langle X,Y\rangle  \de(a+b). $$
Here $\langle \,,\,\rangle $ is the pairing  (\ref{pair}) in the algebra $gl(V_R)$ and a discrete $\de $-function $\delta(a)$ is defined
in the standard way:
$$
 \de(a)=\left\{\!\!
\begin{array}{ccl}
1 &{\rm if}& a=0\\
\rule{0pt}{5mm}
0&{\rm if}&  a\not=0.\\
\end{array}
\right.
$$

We do not know what is the Jacobi identity in  the braided Lie algebra $\agl $, provided $R$ is a Hecke symmetry. However, if $R$ is
an involutive symmetry, the corresponding Jacobi identity is similar to (\ref{Jac}). This claim can be easily deduced from the
following property of the term $\omega$.

\begin{proposition} If $R$ is an involutive symmetry, then the following holds
$$
\om\,[\,,\,]_{23} \left( (I+\Ren_1 \Ren_2+\Ren_2 \Ren_1) (X[a]\ot  Y[b]\ot Z[c])\right)=0.
$$
\end{proposition}

In virtue of this property the term $\omega$ can be called {\em braided cocycle}.

The enveloping algebra  $U(\agl)$ can be also defined in a natural way as the quotient of the free tensor algebra of $\agl $ over the ideal,
generated by the elements
$$
c\, l_i^j[a]-l_i^j[a]\, c,\quad X[a] Y[b]-\RRR(X[a]\ot Y[b])-[X[a] ,Y[b]]-\om(X[a], Y[b])c.
$$
In a similar manner it is possible to define the enveloping algebra  $U(sl(V_R))$. In fact, we suggest a new way of introducing quantum
analogs of affine Lie algebras.

In our subsequent publications we plan to study the center of the  algebras  $U(gl(V_R))$ and $U(sl(V_R))$ in the spirit of the Kac's approach.

\section{Characteristic polynomials for generating matrices}

In this section we suppose that the bi-rank of a given skew-invertible Hecke symmetry $R$ is $(m|0)$, $m\geq 2$. As was shown in
\cite{GPS1}, the generating matrix $L$ of the corresponding RE algebra meets the following Cayley-Hamilton identity
\be
L^m-q L^{m-1} e_1(L)+q^2 L^{m-2}e_2(L)+...+(-q)^{m-1} L e_{m-1}(L)+(-q)^{m} I e_m(L)=0,
\label{CH}
\ee
where
$$
e_0(L)=1,\qquad  e_k(L):=Tr_{R(1\dots k)} (\P_-^{(k)} L_{\overline 1} \,L_{\overline 2}\dots L_{\overline k}),\quad k\geq 1,
$$
are quantum analogs of the elementary symmetric polynomials. Here $\P_-^{(k)}: V^{\ot k} \to V^{\ot k}$ is the skew-symmetrizer (\ref{skew}).

Note that quantum analogs of these and other symmetric polynomials (Schur polynomials, power sums) are also well defined in all QM
algebras and they generate a commutative subalgebra called {\em characteristic}. As we said above, in the RE algebra the characteristic subalgebra
is central. By this reason, it does not matter on what side of the powers of the matrix $L$ we put the coefficients in the Cayley-Hamilton identity
and in  its generalizations called Cayley-Hamilton-Newton identities. An important consequence of this fact is the possibility to
introduce the quantum spectrum of $L$, the quantum eigenvalues $\mu_i$ of $L$ belong to an algebraic extension of the center
of RE algebra. This quantities allows one to rewrite the Cayley-Hamilton identity (\ref{CH}) in a factorized form:
$$
\prod_{i=1}^N(L-\mu_iI) = 0.
$$

By contrary, in other QM algebras  the elements $e_k$ are not central and their position in the Cayley-Hamilton identity (in front of
a matrix power or behind it) is important. Besides, in these algebras the powers $L^k$ should be replaced by their quantum counterparts,
which also exist in two forms.  We refer the reader to the paper \cite{IOP}, where the Cayley-Hamilton identity (and its generalization called the
Cayley-Hamilton-Newton identity) is proved for the generating matrices of all QM algebras, the RTT and RE algebras included.

\begin{remark} \rm In \cite{IO} the notion of a half-quantum algebra was introduced. Similarly to a QM algebra a half-quantum algebra is
defined  with the help of a couple $(R,F)$ of compatible braidings, but the defining relations on generators are less restrictive than those
for the QM algebras. The point is that analogs of symmetric polynomials can be also defined in half-quantum algebras and a version of the
Cayley-Hamilton-Newton relations can be established. Nevertheless, in general, the "symmetric polynomials" in these algebras do not
commute with each other.
\end{remark}

As follows from formula (\ref{CH}), the characteristic polynomial for the generating matrix  $L$ of the RE algebra is
\be
Q(t)=\sum_{k=0}^m t^{m-k} (-q)^k e_k(L),
\label{har-pol}
\ee
since $Q(t)$ is the $m$-th order polynomial with the unit coefficient at the highest power $t^m$ which possesses the property
$Q(L)\equiv 0$.


Our current aim is to present this polynomial in a form useful for finding the characteristic polynomial for the generating matrix of
the modified RE algebra.  By passing to the limit $q\to 1$, we get a characteristic polynomial for the generating matrix  of the algebra $U(gl(N))$.

\begin{proposition}
The polynomial $Q(t)$ defined in {\rm (\ref{har-pol})} is identically equal to the expression:\rm
\be
Q(t)=q^m \,\Tr_{R(1...m)} \left(\P_-^{(m)}(tI-L_{\overline 1})(q^2tI-L_{\overline 2})...(q^{2(m-1)} tI-L_{\overline m})\right).
\label{chp}
\ee
\em
\end{proposition}

\noindent
{\bf Proof.} Consider the following polynomial in $m$ indeterminates $t_i$:
$$
{\hat Q}(t_1,\dots,t_m) =q^m\, \Tr_{R(1\dots m)}\left(\P_-^{(m)}(t_1I-L_{\overline 1})(t_2I-L_{\overline 2})\dots (t_mI-L_{\overline m})\right).
$$
By developing the product of linear factors in the above expression, we get a sum with the typical term
$$
q^m\sigma_k(t_1,\dots ,t_m)\, \Tr_{R(1\dots m)}\left(\P_-^{(m)} (-L_{\overline 1})\dots (-L_{\overline{m-k}})\right),
$$
where $\sigma_k(t_1,\dots ,t_m)$ are the elementary symmetric polynomials in $t_1,...,t_m$:
$$
\sigma_k(t_1,\dots ,t_m) = \sum_{1\le i_1<\dots <i_k\le m}t_{i_1}\dots t_{i_k}.
$$
Here, we use an essential fact that the factors $L_{\overline i}$ {\it under the $R$-trace} can be ``shifted'' to the most possible left position.
 This means that the following identities hold
$$
\Tr_{R(1\dots m)}\left(\P_-^{(m)}L_{\overline s_1}L_{\overline s_2}\dots L_{\overline s_k}\right)=
\Tr_{R(1\dots m)}\left(\P_-^{(m)}L_{\overline 1}L_{\overline 2}\dots L_{\overline k}\right)
$$
for any set of integers $1\leq s_1<s_2 <\dots <s_k\leq m$. Emphasize that this property is specific for the RE algebra. It fails if the
matrices $L_{\overline k}$ are defined via a braiding $F$ different from $R$ (see footnote \ref{fo}).

By using Proposition 1, we can present  the polynomial ${\hat Q}$ as follows
\be
{\hat Q}(t_1,\dots, t_m)=\sum_{k=0}^m  (-1)^k q^{-m(m-k-1)} \,\frac{k_q! (m-k)_q!}{m_q!} \,\sigma_{m-k}(t_1,\dots, t_m) \,e_k(L).
\label{prom-pol}
\ee
Now, we take into account a result from $q$-combinatorics:
$$
\sigma_k(t,q^2t,..., q^{2(k-1)}t)= t^k\sigma_k(1,q^2,..., q^{2(k-1)}) =q^{k(m-1)} \frac{m_q!}{k_q! (m-k)_q!}\,t^k.
$$
So, specializing $t_k = q^{2(k-1)}t$ in the above expression (\ref{prom-pol}) we precisely get the formula (\ref{har-pol}). Therefore
$$
Q(t)={\hat Q}(t,q^2t,\dots, q^{2(m-1)}t),
$$
which proves the claim in virtue of the the definition of $\hat Q$.\hfill \rule{6.5pt}{6.5pt}

\medskip

It is obvious, that formula (\ref{chp}) can be written as follows:
$$
Q(t)=q^m\,\Tr_{R(1\dots m)} \left(\P_-^{(m)}(tI-L)_{\overline 1}(q^2tI-L)_{\overline 2}\dots (q^{2(m-1)} tI-L)_{\overline m}\right).
$$
From this form of the characteristic polynomial (\ref{chp}) for the matrix $L$ we can get the characteristic polynomial for the matrix $\LL$.

\begin{corollary}
\label{cor-1}
The characteristic polynomial ${\tilde Q}(t)$ for the matrix $\LL $ is equal to\rm
\be
{\tilde Q}(t)=q^m\, \Tr_{R(1\dots m)} \left(\P_-^{(m)}\prod_{k=1}^m\left(q^{2(k-1)}\Big(t - q^{-k+1}(k-1)_q\Big)I-\LL_{\overline k}\right)\right),
\label{chp1}
\ee
\em
where the factors are placed in ascending order in $k$ from the left to right.
\end{corollary}

\noindent
{\bf Proof.} The generating matrices $L$ and $\LL $ are connencted by a linear shift (\ref{iso}).  Introducing a new indeterminate
$$
{\tilde t}=t+\frac{1}{q-\qq},
$$
we obviously have $tI - L = \tilde tI -\LL $. Therefore, rewriting $L$ and $t$ in (\ref{chp}) in terms of $\LL $ and $\tilde t$ we get
the polynomial $\tilde Q(\tilde t)$ with the property $\tilde Q(\LL)= 0$  (we return to the notation $t$ instead of $\tilde t$ at the end of transformations).\hfill \rule{6.5pt}{6.5pt}

\medskip

Note, that the value $m$ in the above formulae is in general independent of the parameter $N=\dim V$ except
for the  restriction $m\le N$. If a given  Hecke symmetry  is a deformation of the usual flip
(the $U_q(sl(N))$ Drinfrld-Jimbo $R$-matrix as a well-known example), then its bi-rank is $(N|0)$. By passing to the  limit $q\to 1$ we get the following claim.

\begin{corollary}
Let $M=\|m_i^j\|$ be the generating matrix of the enveloping algebra $U(gl(N))$, where $m_i^j$ $1\le i,j\le N$ is the standard
basis of $gl(N)$. Then the following polynomial is characteristic for this matrix
$$
{\cal Q}(t)= \Tr_{(1\dots N)} \left(\P_-^{(N)}(tI-M_1)((t-1) I-M_2)...((t-N+1) I-M_N)\right),
$$
where $\P_-^{(N)}$ is the usual skew-symmetrizer in $V^{\otimes N}$, that is we have ${\cal Q}(M)=0$.
\end{corollary}

Emphasize that usually the characteristic polynomial for the matrix $M$ is constructed by means of the Capelli determinant.

Also, note that the claim of the corollary \ref{cor-1} is valid for the generating matrix of any modified RE algebra corresponding to a
skew-invertible involutive symmetry $R$, provided $R$ can be approximated by a family of Hecke symmetries.

\begin{remark}\rm
Along with the characteristic polynomial for the generating matrix of the enveloping algebra $U(gl(N))$ one usually exhibits a similar polynomial
for its transposed matrix $M^t$ (see \cite{Mo}). In our setting the matrix $M^t$ should be replaced by the generating  matrix  of the modified
RE algebra defined by
\be
R_1 L_2 R_1L_2-L_2 R_1 L_2 R_1=R_1L_2-L_2 R_1.
\label{mRE1}
\ee
It is treated as the enveloping algebra of the algebra of right endomorphisms of  $V$. For the generating matrix of
this algebra the characteristic  polynomial is similar to (\ref{chp}) but the $R$-trace should be defined via
$\Tr_R L=\Tr (B L)$ and all the matrices $L_{\overline k}$ should be replaced by $L_{\underline k}$ where
 $$
 L_{\underline m}=L_m, \quad L_{\underline{k}}=R_{k}^{-1} L_{\underline {k+1}}R_{k}, \,\, 1\le k\le m-1.
 $$
 \end{remark}

\section{Braided Yangians}
Let $R(u,v)$ be a current quantum $R$-matrix. This means that it is subject to the
quantum Yang-Baxter equation
$$
R_{12}(u,v) R_{23}(u, w) R_{12}(v,w)  =R_{23}(v,w) R_{12}(u, w) R_{23}(u,v).
$$
Consider an analog of the RTT algebra, associated with such an $R$-matrix defined by the system
\be
R(u,v) L_1(u)  L_2(v)=L_1(v)  L_2(u)R(u,v),
\label{rtt}
\ee
where the matrix
\be
L(u)=\sum_{k\geq 0} L[k] u^{-k},\qquad L[k]=\|l_i^j[k]\|_{1\leq i, j \leq N},
\label{expa}
\ee
expands in a series in non-positive powers of the parameter.

Thus, the system (\ref{rtt}), being rewritten via the generators $l_i^j[k]$, leads to an infinite family of quadratic equations on the
generators $l_i^j[k]$, whose number is also infinite, but each equation is a polynomial in the generators.

The $R$-matrix  $R(u,v)$, we are dealing with, are of two classes:
\be
R(u,v)=R-\frac{a}{u-v}I,\qquad R(u,v)=R-\frac{(q-\qq)u}{u-v}.
\label{Rm}
\ee

In the first formula in (\ref{Rm}) $R$ stands for a skew-invertible involutive symmetry, whereas in the latter one $R=R(q)$ is a
skew-invertible Hecke symmetry. The fact that these $R$-matrices meet the quantum Yang-Baxter equation can be verified by a
straightforward calculation \cite{GS1}. The procedure of constructing such current $R$-matrices via symmetries $R$ is often called
{\em Yang-Baxterization}.

Note that the Drinfeld's Yangian corresponds to the famous Yang $R$-matrix, which is defined via the first formula (\ref{Rm}) with $R=P$.

Observe that the second current $R$-matrix in (\ref{Rm}) is actually depending on the ratio  $u/v$ of the parameters.
There exists another (so-called {\em trigonometrical}) form which can be obtained via the change of variables $u\mapsto q^u$, $v \mapsto q^v$.
After such a transformation the current $R$-matrix will depend on the difference $u-v$. Below, we deal with the rational form (\ref{Rm}) of the
current $R$-matrix.

All Yangian-like algebras defined via (\ref{rtt}) but with other current $R$-matrices are called {\em Yangians of RTT type}.

The well-known examples of  Yangians of RTT type are the so-called $q$-Yangians (see \cite{Mo}). In the lowest dimensional case such a $q$-Yangian is associated with
the following current $R$-matrix
\be
R(u,v)= R-\frac{(q-\qq)u}{u-v}I=
\left(\begin{array}{cccc}
\frac{-qv+\qq u}{u-v}&0&0&0\\
0&\frac{(-q+\qq)v}{u-v}&1&0\\
0&1&\frac{(-q+\qq)u}{u-v}&0\\
0&0&0&\frac{-qv+\qq u}{u-v}
\end{array}\right).
\label{Rm1}
\ee

Note that each Yangian of RTT type has a bi-algebra structure. On the generators the corresponding coproduct is defined as follows
$$
\De(1)=1\ot 1,\qquad \De(l_i^j(u))=l_i^k(u)\ot l_k^j(u).
$$

Now, consider the so-called {\em evaluation morphism}
$$
T(u)\to T+\frac{\OT}{u}\,.
$$
This map induces a morphism of algebras, if $R$ is a Hecke symmetry and the matrices $T$ and $\OT $ meet the following relations
$$
R\,T_1\, T_2=T_1\, T_2\, R,\quad R\,\OT_1\, \OT_2=\OT_1\, \OT_2\, R,\quad R\,\OT_1\, T_2=T_1\, \OT_2\, R.
$$
Thus, the target algebra generated by the entries of the matrices $T$ and  $\OT $ is a couple of RTT algebras connected by the last relation.

Note that in the case of $q$-Yangians one usually imposes some additional conditions on the matrix $L[0]$.

In \cite{GS1} we suggested another candidate for the role of the $q$-Yangian. Consider an algebra generated by entries of a matrix $L(u)$
subject to the relations
\be
R(u,v) L_1(u) R L_1(v)-L_1(v) R L_1(u) R(u,v)=0,
\label{brY}
\ee
where $R$ is just the involutive or Hecke symmetry entering the current $R$-matrix (\ref{Rm}). Besides, in the expansion  (\ref{expa}) we assume that $L[0]=I$.
We denote this algebra $\YR$ and call it the {\em braided Yangian}.

Below, we use the notation similar to (\ref{ove}):
$$
L_{\overline 1 }(u)=L_{1}(u),\quad L_{\overline k}(u)=R_{k-1} L_{\overline{k-1}}(u)R_{k-1}^{-1},\quad k\ge 2.
$$
By using this notation, it is possible to cast the defining relations (\ref{brY}) in a form similar to the Yangians of RTT type
$$
R(u,v) L_{\overline 1}(u) L_{\overline 2}(v)-L_{\overline 1}(v)  L_{\overline 2}(u)  R(u,v)=0.
$$

The braided Yangian, corresponding to the $R$ matrix (\ref{Rm1}) and its higher dimensional analogs are called {\em braided $q$-Yangian}.

Let us mention some of properties of the braided Yangians. First of all, any braided Yangian has a braided bi-algebra structure.
The corresponding coproduct is defined on the generators by the same formulae as in the Yangians of RTT type but it is extended
on the whole algebra via the braiding $\Ren $ in a way similar to that in the RE algebra.

Another important property of the braided Yangians is that their evaluation morphisms look like the classical one
\be L(u) \to I+\frac{M}{u}. \label{ev} \ee
The target algebra generated by entries of the matrix $M$, is described by the following claim proved in \cite{GS1}.

\begin{proposition}

\hspace*{5mm}

\begin{enumerate}
\item[\rm 1.] If $R$ is an involutive symmetry, then the map {\rm (\ref{ev})} defines a surjective morphism $\YR\to \tlr $. Besides, the map
$M\mapsto L[1]$ defines an injective morphism $\tlr\to \YR$.

\item[\rm 2.] If $R$ is a Hecke symmetry, then the map {\rm (\ref{ev})} defines a morphism $\YR\to \lr $.
\end{enumerate}
\end{proposition}

Thus, the type of the target algebra depends on the type of the initial symmetry $R$. This proposition enables us to construct a large
representation category of each braided Yangian. We describe it for the braided $q$-Yangian $\YR $.

Consider the category of finite dimensional $U_q(sl(N))$-modules which are deformations of the $U(sl(N))$-ones. Each of its
objets can be endowed with a structure of the $\tlr $-module where $R$ is coming from the QG  $U_q(sl(N))$. This fact follows
from the method of constructing the category of $\tlr $-module as was done in \cite{GPS1}. Finally, by using the isomorphism (\ref{iso})
we can convert any  $\tlr $-module into $\lr $-one. Now, it remains to apply the above proposition.

One of the most remarkable properties of the Yangians of all types is that analogs of some symmetric functions are well defined in these algebras.
Also, analogs of the Cayley-Hamilton-Newton identities are valid in these algebras. Note that these identities can be presented in different
form. Below, we exhibit them in a form which differs from that of \cite{GS1,GS2}.

First, define analogs of  powers $L^k(u)$ and skew-powers $L^{\wedge k}(u)$ of the matrix $L(u)$, generating a braided Yangian:
$$
\begin{array}{l}
L^k(u):= L(q^{-2(k-1)} u) L(q^{-2(k-2)} u)\dots L(u),\\
\rule{0pt}{5mm}
L^{\wedge k}(u): = \Tr_{R(2\dots k)}\Big({\cal P}^{(k)}_{-}L_{\overline 1}(u)L_{\overline 2}(q^{-2} u)\dots L_{\overline k}(q^{-2(k-1)}u)\Big)
\quad k\ge 2,
\end{array}
$$
where it is convenient to set by definition $L^0(u) = I$ and $L^{\wedge 1}(u) = L(u)$.

Here, we would like to emphasize a difference between the braided Yangians and these of RTT type. In the latter ones analogs of the
matrix powers and skew-powers can be also defined, but only in the braided Yangians analogs of the matrix powers are defined via the usual
matrix product of the generating matrices (but with shifted parameters).

In the braided Yangians quantum analogs of the power sums and elementary symmetric polynomials are respectively defined by
\begin{eqnarray}
p_k(u)&=&\Tr_R L^k(u)=\Tr_R L(q^{-2(k-1)} u) L(q^{-2(k-2)} u)\dots L(u),\nonumber\\
e_k(u)&=&\Tr_R L^{\wedge k}(u)=\Tr_{R(1\dots k)}\Big({\cal P}^{(k)}_{-}L_{\overline 1}(u)L_{\overline 2}(q^{-2} u)\dots
L_{\overline k}(q^{-2(k-1)}u)\Big).
\label{elsym}
\end{eqnarray}
Now, we are able to exhibit the Cayley-Hamilton-Newton identities in the braided Yangians.

\begin{proposition}
The following matrix identities hold true for the generating matrix of a braided Yangian\rm
\be
(-1)^{k+1} k_q L^{\wedge k}(u)=\sum_{p=1}^k (-q)^{k-p} L^p(q^{2(k-p)}u) e_{k-p}(u)\quad \forall\,k\ge 1.
\label{chn-id}
\ee\it
\end{proposition}

\medskip

\noindent
{\bf Proof}  can be done by the considerations similar to these from \cite{GS1}, namely, by sequential use of
the recurrent formula (\ref{skew}). However, first, we  present
 the skew-power $L^{\wedge k}(u)$ as follows
$$
L^{\wedge k}(u) = \Tr_{R(2\dots k)}\Big(L_{\overline 1}(q^{-2(k-1)}u)L_{\overline 2}(q^{-2(k-2)}u)\dots L_{\overline k}(u)
{\cal P}^{(k)}_{-}\Big).
$$
Then we apply the method used in \cite{GS1}.

Note, that in the cited paper we got a different form of the
Cayley-Hamilton-Newton identities: the elementary symmetric functions appeared there on the {\it left} of matrix
powers of $L$ and the matrix powers were defined by more complicated expressions than those written above.
\hfill\rule{6.5pt}{6.5pt}
\medskip

If the bi-rank of the Hecke symmetry $R$ is $(m|0)$ (in particular, for the braided $q$-Yangian it is $(N|0)$, i.e. $m=N$), then the highest
nonzero skew-power is
$$
L^{\wedge m}(u) =q^m e_m(u)\, I,
$$
where $e_m(u)$ is the highest nonzero elementary symmetric polynomial, which is a quantum analog of the determinant.

Consequently, on setting in (\ref{chn-id}) $k=m$ we get the Cayley-Hamilton identity
$$
\sum_{p=0}^m (-q)^p L^{m-p}(q^{-2p}u) e_p(u)=0.
$$

By applying the $R$-trace to the identities (\ref{chn-id}), we get a family of the quantum Newton identities
\begin{eqnarray*}
p_k(u)-q p_{k-1}(q^{-2}u) e_1(u)\!\!\!\!&+&\!\!\!\!(-q)^2 p_{k-2}(q^{-4}u)e_2(u)+\dots\\
&+&\!\!\!\!(-q)^{k-1} p_1(q^{-2(k-1)}u) e_{k-1}(u)+(-1)^k k_q e_k(u)=0\quad \forall\,k\ge 1.
\end{eqnarray*}

Emphasize that the quantum analog of the determinant $e_m(u)$ is central in the braided Yangian $\YR $ for any symmetry $R$, whereas
in the Yangians of RTT type its centrality depends on $R$ (see \cite{GS1}).  The other analogs of elementary symmetric  polynomials
$e_k(u)$, $1\le k\le m$ and power sums  $p_k(u)$, $k\ge 1$, commute with each other and generate a commutative {\it Bethe subalgebra}.
This is true also for the Yangians of RTT type.

\section{Shifted braided $q$-Yangian and its $q=1$ limit}

Let $r(u,v)$ be a classical $r$-matrix, i.e. $gl(N)$-valued  function in parameters $u$ and $v$ (assumed to be rational), which meets the classical Yang-Baxter equation
$$ [r_{12}(u,v),r_{13}(u,w)]+ [r_{12}(u,v),r_{23}(v,w)]+[r_{13}(u,w),r_{23}(v,w)]=0.$$
Suppose  that the map
\be L(u)\ot L(v)\mapsto \{L_1(u),L_2(v)\}=[r(u,v), L_1(u)+L_2(v)], \label{Pb} \ee
is skew-symmetric and consequently it
defines a Poisson bracket. Also, suppose that the matrix-function $L(u)=\|l_i^j(u)\|$ expands in a series (\ref{expa}). Thus, this Poisson bracket,
defined on the commutative algebra ${\rm Sym}(gl(N)[t^{-1}])$, can be expressed via the coefficients $L[k],\, k\geq 0$.

By using the standard $R$-matrix technique,  it is easy to show that the elements $\Tr L^k(u)$ commute with each other with respect to this Poisson bracket:
 \be \{\Tr L^k(u), \Tr L^l(v)\}=0,\quad \forall\, k,l \geq 0 \label{comm} \ee
 for any  values $u$ and $v$.

 The simplest non-constant example  of an  $r$-matrix is the following one
\be r(u,v)=\frac{P}{u-v}. \label{rma} \ee
The corresponding Poisson bracket is an important ingredient of the rational Gaudin model and its quantization.

Let us emphasize  that the elements $\Tr L^k(u)$ (or their $R$-counterparts)  do not commute any more in the enveloping algebra $U(\GG)$ of the Lie algebra $\GG$ defined by  formula
(\ref{Pb}) with the same $r$-matrix (\ref{rma}).  D.Talalaev \cite{T} succeeded in finding a Bethe subalgebra in this enveloping algebra.

Now, consider the map
\be L(u)\mapsto \sum_{k=1}^K \frac{M_k}{u-u_k}, \label{LM} \ee
where  $u_k\in \K$ are $K$ fixed points,
and  the matrices $M_k$
generate $K$ copies of the algebra $U(gl(N))$. This means that each matrix $M_k$ generates the enveloping algebras $gl(N)$ and entries of any two matrices of this family commute with each other.
Then, the map (\ref{LM}) defines a Lie algebra morphism $\GG\to \gg=gl(N)^{\oplus K}$ and consequently a morphism of the enveloping algebras of these Lie algebras.

In fact, the map (\ref{LM}) can be treated as an analog of the evaluation morphism $\Y(gl(N))\to U(gl(N))$, combined with the morphisms $u\mapsto u-u_k$ and the usual coproduct.

The image of the Bethe subalgebra in the algebra $U(\GG)$ under the map (\ref{LM}) is a Bethe subalgebra in the algebra $U(\gg)$.
Some quadratic elements of the latter algebra play the role of the Hamiltonians of the rational Gaudin model. Talalaev's result gives rise to higher
Hamiltonians of this model.

Let us emphasize that constructing a Bethe subalgebra in the algebra $U(\GG)$ was performed by Talalaev via using a Bethe subalgebra in the Yangian $Y(gl(N))$.
It is tempting to replace the $r$-matrix (\ref{rma}) by that corresponding to $R$-matrix (\ref{Rm1}) and upon using  the same method, to find a Bethe subalgebra in
the enveloping algebra of the  corresponding Lie algebra.

Unfortunately, this method  fails though a Bethe subalgebra in the $q$-Yangian of RTT type exists and is known (see \cite{GS1}).
This failure is due to the fact that  in the $q$-Yangian of RTT type it is not possible to find  elements of this  Bethe subalgebra which have a controllable
expansion in the deformation parameter $h$.

We claim that the Talalev's method is still valid in the braided $q$-Yangian. Namely, below we exhibit elements of the Bethe subalgebra of this generalized
Yangian which has the necessary expansion property. This enables us to find a Bethe subalgebra in the enveloping algebra $U(\GG_{trig})$
of the Lie algebra $\GG_{trig}$, which arises (similarly to the rational case) from linear Poisson structure corresponding to this braided $q$-Yangian (see \cite{GS1}).
Besides, an analog of the map (\ref{LM}) can be also found  and used for constructing a Bethe subalgebra in the algebra  $U(\gg)$.

As  claimed in \cite{GS1, GSS}, the quantum     elementary symmetric polynomials $e_k(u)$ defined by (\ref{elsym}) commute with each other in the braided Yangians $\YR$.
Also, we consider the following elements
$$ \he_k(u)= \Tr_{R(1...m)}\Big({\cal P}^{(m)}_{-}L_{\overline 1}(u)L_{\overline 2}(q^{-2} u)\dots L_{\overline k}(q^{-2(k-1)}u)\Big), $$
which differ from  $e_k(u)$ by a modification of the skew-symmetrizers  and the positions where the $R$-traces are applied.
According to Proposition \ref{pr1}  any element $\he_k(u)$ differs from  that $e_k(u)$ by a non-trivial (for a generic $q$) numerical factor.
Consequently, the elements $\he_k(u)$ also commute with each other.

By using    the relation
$$q^{-2\pa_u }f(u)=f(q^{-2} u)q^{-2\pa_u },\,\,\,{\rm where}\,\,\,   \pa_u=u\frac{d}{du},$$
we can present the elements $\he_k( q^{-2}u)$ as follows
\be \he_k( q^{-2} u)= \Tr_{R(1...m)}\Big({\cal P}^{(m)}_{-}(q^{-2\pa_u }{L}_{\overline 1}(u))(q^{-2\pa_u }{L}_{\overline 2}( u))\dots  (q^{-2\pa_u }{L}_{\overline k}(u))\Big)q^{2k \pa_u }. \label{hek} \ee

Now,  change  the basis of the braided $q$-Yangian in a way similar to  (\ref{iso}):
\be L(u)=(q-\qq) \tLL(u)+I. \label{chan} \ee
Also, we put $q^2=\exp(h)$ and expand the elements $\he(q^{-2} u)$, expressed via the matrix $\tLL$, in  $h$.

Since the factors entering formula (\ref{hek}) expand as
$$q^{-2\pa_u }{L}_{\overline p}(u)=(1-h\pa_u+o(h))(I+h\tLL_{\overline p}(u)+o(h))= I+h(\tLL_{\overline p}(u)-I\,\pa_u)+o(h),$$
we get the following expansion of the elements $\he_k( q^{-2} u)$
$$\he_k(e^{-h} u)= \Tr_{R(1...m)}\Big({\cal P}^{(m)}_{-}(I+h(\tLL_{\overline 1}(u)-I\,\pa_u)+o(h))...(I+h(\tLL_{\overline k}(u)-I\,\pa_u)+o(h))\Big)(1+o(1)). $$

Also, note that $\tLL_{\overline p}(u)=\tLL_p(u)+o(1)$ for all $p$.

Following \cite{T}, consider the elements
$$\tau_k(u)=\sum_{p=0}^k(-1)^{k-p}\frac{k!}{p! (k-p)!}\he_p(u),$$
commuting with each  other in the braided $q$-Yangian.

We state that the expansion of the element $\tau_k(e^{-h}u)$  begins with a term proportional to $h^{k}$.
Thus, the elements $h^{-k}\tau_k(e^{-h}  u)$ expand as follows
$$h^{-k}\tau_k(e^{-h} u)=QH_k(u)+o(h).$$

This entails that the elements  $QH_k(u), k=0,1...,m$    commute with each other in the algebra which is the limit of the braided $q$-Yangian as $h\to 0$.
Let us compute the defining relations of the limit algebra. Being expressed via the matrix
$\tLL$ the defining system of the braided $q$-Yangian is
 $$(R-\frac{(q-\qq)u}{u-v} I) {\tLL(u)}_1 R {\tLL(v)}_1-{\tLL(v)}_1 R {\tLL(v)}_1 (R-\frac{(q-\qq)u}{u-v} I)=
 -[R, \frac{u\tLL_1(u)-v\tLL_1(v)}{u-v}].$$

Consequently,  the defining relations of the limit algebra are
\be [\tLL_1(u),  \tLL_2(v)]=[P, \frac{u\tLL_1(u)-v\tLL_1(v)}{u-v}]=[\frac{P}{u-v}, u\tLL_1(u)+v\tLL_2(v)]. \label{PS} \ee

We denote $\GG_{trig}$ the  Lie algebra with the bracket, defined by the right hand side of this formula. Thus, in the basis $\tLL(u)$  the braided $q$-Yangian turns  into
 the enveloping algebra $U(\GG_{trig})$ of this  Lie algebra as $q\to 1$ (or $h\to 0$).

Thus, similarly to \cite{T} we have the following.

\begin{proposition} The elements
$$QH_k(u)=\Tr_{(1...m)}\P^{(m)}_- (\tLL_1(u)-I\pa_u)(\tLL_2(u)-I\pa_u)...(\tLL_k(u)-I\pa_u) 1, \,\,\, 1\geq k \geq m$$
commute with each other in the algebra $U(\GG_{trig})$.   \end{proposition}

Observe that here the trace and skew-symmetrizer are classical.

The commutative subalgebra generated by the elements $QH_k(u)$ in the algebra $U(\GG_{trig})$ is called the {\em Bethe subalgebra}.

\begin{remark} \rm Note that similarly to (\ref{iso})
the  map (\ref{chan}) converts a quadratic algebra in a quadra\-tic-linear one. However,  the defining relations of the limit algebra
do not depend on
the concrete form of the Hecke symmetry $R$.  \end{remark}

Now, describe an analog of the map (\ref{LM}).

\begin{proposition} The map
$$ \tLL(u)\mapsto \sum_{k=1}^K \frac{ M_k u_k}{ u-u_k},$$
where the family $(M_1,..., M_K)$ generates the Lie algebra $\gg=gl(m)^{\oplus K}$ (in the same sense as in (\ref{LM})),
defines a Lie algebra morphism $\GG_{trig} \to \gg$ and consequently, a morphism of the enveloping algebras of these Lie algebras.
\end{proposition}

The corresponding version of an integrable model will be considered elsewhere.

Let us point out two main differences of our result from that of \cite{T}.
First,  our  $\pa_u$ is the "multiplicative derivative": $u\frac{d}{du}$. Second, the Lie algebra $\GG_{trig}$ cannot be presented in the form (\ref{Pb}).
The Poisson structure, corresponding to this Poisson bracket, is exhibited in  \cite{GS2}. Also,  a comparative analysis of the Poisson structures related to different types of the deformation Yangians is presented there. Emphasize that the Poisson structures corresponding to the braided Yangians do not depend on the concrete form of the Hecke symmetry
 as well.  Up to a factor, it equals  the right hand side of formula (\ref{PS}).

Completing the paper, emphasize once more that the Talalaev's method is still valid since we are dealing with the braided version of the  Yangians.


\begin{thebibliography}{MRZZ}

\bibitem[D]{D} V.Drinfeld, {\em Quantum Groups}, Proceedings of the International Congress of Mathematicians,
Vol. 1, 2 (Berkeley, Calif., 1986), 798--820, Amer. Math. Soc., Providence, RI, 1987.

\bibitem[FJMR]{FJMR} L.Frappat, N.Jing, A.Molev, E.Ragoucy, {\em Higher Sugawara operators for the quantum affine algebras of type $A$}, arXiv:1505.03667.

\bibitem[G]{G} Gurevich D. {\em Algebraic aspects of quantum Yang-Baxter equation}, Leningrad Math. Journal 2:4 (1990), 119--148.

\bibitem[GPS1]{GPS1} D. Gurevich, P. Pyatov, P. Saponov {\em Hecke symmetries and characteristic relations on reflection equation algebras},  Lett. Math. Phys. 41 (1997), no. 3, 255–264.

\bibitem[GPS2]{GPS2} D. Gurevich, P. Pyatov, P. Saponov {\em Representation theory of (modified) Reflection Equation Algebra of the $GL(m|n)$ type},
 Algebra and Analysis, {\bf 20} (2008), 70--133.

\bibitem[GS1]{GS1} D. Gurevich, P. Saponov {\em Braided Yangians}, arXiv:1612.05929.

\bibitem[GS2]{GS2} D. Gurevich, P. Saponov, {\em Generalized  Yangians and their Poisson counterparts}, arXiv:1702.03223.

\bibitem[GSS]{GSS} D. Gurevich, P. Saponov, A. Slinkin {\em Bethe subalgebras and matrix identities in Braided Yangians}, in progress.

\bibitem[H]{H} Phung Ho Hai, {\em Poincar\'e Series of Quantum spaces Associated to Hecke Operators}, Acta Math. Vietnam  24  (1999) 235–-246.

\bibitem[IP]{IP} A.Isaev, P.Pyatov, {\em Spectral extension of the quantum group cotangent bundle}, Comm. Math. Phys. 288 (2009), no. 3, 1137–-1179.

\bibitem[IO]{IO} A.Isaev, O.Ogievetsky {\em Half-quantum linear algebra} in:  Symmetries and groups in contemporary physics, pp 479–486,
Nankai Ser. Pure Appl. Math. Theoret. Phys., 11, World Sci. Publ., Hackensack, NJ, 2013.

\bibitem[IOP]{IOP} A.Isaev, O.Ogievetsky, P.Pyatov, {\em On quantum matrix algebras satisfying the Cayley-Hamilton-Newton identities}, J. Phys. A 32 (1999), no. 9, L115--L121.

\bibitem[M]{M} Sh.Majid, {\em Foundations of quantum group theory}, Cambridge University Press, Cambridge, 1995.

\bibitem[Mo]{Mo} A.Molev, {\em Yangians and classical Lie algebras}, Mathematical Surveys and Monographs, 143. American Mathematical Society, Providence, RI, 2007.

\bibitem[O]{O} O.Ogievetsky, {\em Uses of Quantum Spaces}, 3rd cycle. Bariloche (Argentine), 2000, pp.72, cel-00374419.

\bibitem[RS]{RS} N.Reshetikhin, M.Semenov-Tian-Shansky, {\em Central extensions of quantum current groups},  Lett. Math. Phys. 19 (1990), no. 2, pp. 133–142.

\bibitem[T]{T} D.Talalaev {\em Quantum Gaudin system}, Func. Anal. Appl. 40 (2006), no 1, 73-77.

\end{thebibliography}
\end{document}